\input amstex
\documentstyle{gen-j}

\chardef\ss="19
\def\3{\ss}
\def\<{<\!\!<}

\def\into{\hookrightarrow}
\def\ssarr{\hbox to 30pt{\rightarrowfill}}
\def\sarr{\hbox to 40pt{\rightarrowfill}}
\def\arr{\hbox to 60pt{\rightarrowfill}}
\def\larr{\hbox to 60pt{\leftarrowfill}}
\def\Arr{\hbox to 80pt{\rightarrowfill}}

{}

\def\Ad{\mathop{\text{\rm Ad}}\nolimits}

\def\co{\mathop{\text{\rm co}}\nolimits}
\def\conv{\mathop{\text{\rm conv}}\nolimits}
\def\diag{\mathop{\text{\rm diag}}\nolimits}

%
%

\def\Im{\mathop{\text{\rm Im}}\nolimits}

\def\pr{\mathop{\text{\bf pr}}\nolimits}

\def\Re{\mathop{\text{\rm Re}}\nolimits}
\def\Res{\mathop{\text{\rm Res}}\nolimits}

\def\Sl{\mathop{\text{\rm Sl}}\nolimits}
\def\SU{\mathop{\text{\rm SU}}\nolimits}

\def\sup{\mathop{\text{\rm sup}}\nolimits}

\def\0{{\bf 0}}
\def\1{\text{\bf {1}}}

\def\a{{\frak a}}

\def\c {\mathop{\text{\bf c}}\nolimits}

\def\g{{\frak g}}

\def\k{{\frak k}}
\def\l{{\frak l}}
\def\m{{\frak m}}

\def\n{{\frak n}}

\def\p{{\frak p}}

\def\z{{\frak z}}

\def\C{{\Bbb C}}

\def\K{{\Bbb K}}

\def\N{{\Bbb N}}

\def\R{{\Bbb R}}

\def\:{\colon}  
\def\.{{\cdot}}
\def\|{\Vert}

\def\giantskip{\vskip2\bigskipamount}

\def \la {\langle}

\def \ra {\rangle}
\def \res {\!\mid\!\!}

\def\giantbreak{\par \ifdim\lastskip<2\bigskipamount \removelastskip
         \penalty-400 \giantskip\fi}

\def\nin{\noindent}
\def\pagebreak{\vskip 0pt plus 0.0001fil\break}
\def\linebreak{\break}

\def\hat{\widehat}

\def\nin{\noindent}
\def\oline{\overline}

\def\phi{\varphi}


\def\subeq{\subseteq}

\def\tilde{\widetilde}

\def\into{\hookrightarrow}
\def\onto{\to\mskip-14mu\to} 

\def\bs{\backslash}

\def\SO{\mathop{\text{\rm SO}}\nolimits}

\topmatter 
\title A convexity property for real orbits in complexified 
Riemannian symmetric spaces \endtitle

\rightheadtext{Convexity in complexified Riemannian symmetric spaces}
\leftheadtext{Bernhard Kr\"otz}

\author Bernhard Kr\"otz \endauthor 

\address The Ohio State University, Department of Mathematics, 
231 West 18th Avenue, Columbus OH 43210-1174\endaddress 
\email kroetz\@math.ohio-state.edu\endemail

\thanks The work was supported in part by NSF-grant DMS-0097314\endthanks

\endtopmatter 
\document 

\head Introduction\endhead 

For real semisimple Lie groups $G$ and their Lie algebras $\g$ interesting 
features of convexity have been discovered. This paper makes 
a contribution which is related to certain $G$-orbits in the 
complexification $X_\C=G_\C/ K_\C$ of a Riemannian symmetric 
space $X=G/K$ of the non-compact type. Applications 
to analytically continued spherical functions will be given.

\par Our results borrow in spirit and techniques from the 
well understood situation on the real space $X$. Therefore it is useful 
to review some of these results and methods. Let $\g$ and $\k$ 
denote the Lie algebras of $G$ and $K$. Then there is a 
Cartan decomposition $\g=\k+\p$. We fix a maximal 
abelian subspace $\a\subeq\p$. Then 
{\it Kostant's linear convexity theorem} [5] reads 
$$(\forall Y\in \a)\qquad  \pr(\Ad(K)Y)=\conv({\Cal W}Y)\ .\leqno(1)$$
In this formula $\pr\: \p\to \a$ is the orthogonal projection, 
${\Cal W}=Z_K(\a)\bs N_K(\a)$ is the Weyl group and $\conv(\cdot)$
denotes the convex hull of $(\cdot)$. 
\par There is a related version of (1) on the Riemannian symmetric 
space $X=G/K$. For that let $G=NAK$ be an Iwasawa decomposition of $G$ and 
$a\: G\to A$ the corresponding middle projection. Then 
{\it Kostant's non-linear convexity theorem} [5] says 
$$(\forall Y\in \a)\qquad  \log a(K\exp(Y))=\conv({\Cal W}Y)\ .\leqno(2)$$
\par  Let us mention that for all applications of (1) and (2) to harmonic
analysis one only needs the inclusion $''\subeq''$. 
For this part (1) and (2) can easily seen to be equivalent: 
Let $\lambda\in \a^*$ and $Y\in \a$ be both regular. Then with $M=Z_K(\a)$ 
the linear/non-linear phase functions on $M\bs K$ are defined by 

$$f_{\lambda, Y}\: M\bs K\to \R, \ \ Mk\mapsto \lambda(\pr(\Ad(k)Y))$$
$$F_{\lambda, Y}\: M\bs K\to \R, \ \ Mk\mapsto a(k\exp(Y))^\lambda\ .$$
Standard arguments show that $f_{\lambda, Y}$ and $F_{\lambda, Y}$
have the same critical set which is computed to be ${\Cal W}$ (cf.\ [12]). 
From that one easily derives both inclusions $``\subeq''$ in (1) and (2).

\par The phase function $F_{\lambda, Y}$  (as well as $f_{\lambda, Y}$)
is Morse. This means that the method of stationary phase applies for computing 
the asypmtotics of the oscillatory integral associated 
to $F_{\lambda, Y}$ (cf.\ [12]). 
Up to a spectral shift by $\rho\in \a^*$ this  oscillatory integral 
is the spherical function $\phi_\lambda$ on $X=G/K$

$$\phi_\lambda(a)=\int_{M\bs K} F_{\rho -i\lambda, Y} (Mk)\ dMk \qquad 
(a=\exp(Y))\ .$$
The conclusion is that convexity associated to the Riemannian 
symmetric space $X=G/K$ is intimately related to the asymptotics 
of spherical functions in the spectral parameter $\lambda$. 
It is important for us that this dictionary remains true once 
we analytically continue from $X$ to certain domains in $X_\C$
(although the method of stationary phase does not apply anymore). 
We will comment more on this later.

\medskip  There is another convexity result for $X$ which 
is related to the polar decomposition $X=K\oline {A^+}K/K$:  
$$(\forall Y_1, Y_2\in \oline{\a^+})\qquad \exp(Y_1)K\exp(Y_2)\subeq K
\exp\left(\conv\left({\Cal W}(Y_1 +Y_2)\right)\right) K\  \leqno(3)$$ 
(cf.\ [11]) 
\par The main result of this paper is a certain complex version of (3). 
Let us be more specific. Denote by $\Sigma$ the root system of $\g$ w.r.t. 
$\a$ and define a bounded convex domain  in $\a$ by 
$$\Omega=\{ Y\in \a\: (\forall \alpha\in \Sigma)\ |\alpha(Y)|<{\pi\over 2}\}
\.$$ 
If $x_o$ denotes the base point in $X_\C$, then we 
shall be interested in those $G$-orbits in $X_\C$ which pass
through $\exp(i\Omega).x_o$. The union of those orbits 
gives us the {\it complex crown} of the Riemannian symmetric space 
$X$:
$$\Xi=G\exp(i\Omega).x_o\subeq X_\C\ .$$ 
The complex crown is a $G$-invariant Stein neighborhood of $X$ in $X_\C$ 
which was subject of many recent investigations. 
\par Finally we set $X_{\C, s}=K_\C A_\C.x_o$ and notice that 
$X_{\C,s}$ is set of all points in $X_\C$ with closed $K_\C$-orbits. 
This set is large in the sense that it contains an open and dense 
subset of $X_\C$ (regular semisimple elements).  
The appropriate holomorphic extension of (3) now is
(cf.\ Theorem 2.1 and Corollary 2.2 below):

\proclaim{Theorem 1} Let $Y\in \Omega$. Then 
$$G\exp(iY).x_o\bigcap X_{\C, s} \subeq K_\C A\exp(i\conv({\Cal W}Y)).x_o\
 .$$
In particular, $\Xi\bigcap X_{\C, s}\subeq K_\C 
A\exp(i\Omega).x_o\ .$\qed \endproclaim 

\medskip For the Lorentzian groups $G=\SO_e(1,n)$ this 
convexity inclusion was proved  in [7] by explicit computation. In the 
classical cases of $G=\Sl(n,\K)$, i.e. when 
$X$ has a canonical realization in symmetric matrices, it might be 
interesting to observe that 
Theorem 1 is not obvious from elementary linear algebra (see Section 3
for a detailed discussion).

\par We describe briefly the method of the proof. 
It was observed in [8] that spherical functions $\phi_\lambda$ 
on $X$ admit a holomorphic continuation $\tilde\phi_\lambda$ to 
$\Xi$. Equally spherical functions extend to $K_\C$-invariant 
holomorphic functions on $K_\C A\exp(2i\Omega).x_o$ (cf.\ [8]).
For fixed generic elements $a\in A\exp(i\Omega)$ we give near optimal 
estimates for $\tilde\phi_\lambda(a)$ for large spectral parameters 
$\lambda$. The upper estimate hereby is obtained from the complex version 
of the non-linear convexity theorem [2]. For the lower 
estimate we use Harish-Chandra's power series expansion for 
spherical functions. From these estimates one then derives the 
convexity inclusion of Theorem 1.   
\par The above convexity theorem features various applications to 
harmonic analysis. In summary we show (cf.\ Theorem 2.3.1 and 
Theorem 2.3.2 below): 

\proclaim{Theorem 2} The following assertions hold: 
\roster 
\item Let $\phi_\lambda$ be a positive definite spherical function 
on $X$. Then its holomorphic extension $\tilde\phi_\lambda$ 
to $\Xi$ is bounded. 
\item Let $\rho_t$ be the heat kernel on $X$ and $\tilde\rho_t$ its 
holomorphic continuation to $\Xi$  (cf.\ [9]). Then $\tilde\rho_t$
is bounded. 
\endroster\qed 
\endproclaim

\medskip The author would like to thank M. Otto, C. Rader and 
M. Stenzel for several useful discussions.

\head 1. Basic facts on $K_\C$-double cosets and complex crowns\endhead 

This section is organized as follows. We begin with introducing  the 
standard notation neded for the discussion 
of Riemannian symmetric spaces $X=G/K$ 
and their complexifications $X_\C=G_\C/K_\C$. Then we summarize 
the basic facts of the $K_\C$
-action on the affine variety $X_\C$. Finally we provide the 
reader with the definition and main features of the complex crown 
$\Xi\subeq X_\C$ of a Riemannian symmetric space $X$.

\subhead 1.1. Notation \endsubhead 
Let us denote by $\g$ a real semisimple Lie algebra with 
Cartan decomposition $\g=\k+\p$. We fix a maximal 
abelian subspace $\a\subeq \p$ and set $\m=\z_\k(\a)$. 
Then one  has the root space decomposition 
$$\g=\a\oplus\m \oplus \bigoplus_{\alpha\in \Sigma} \g^\alpha$$
with $\Sigma\subeq \a^*$ the corresponding restricted root system. 
Select a positive system $\Sigma^+\subeq \Sigma$ and form the 
nilpotent subalgebra $\n=\bigoplus_{\alpha\in \Sigma^+} \g^\alpha$.

\par If $\l$ is a Lie algebra, then we write $\l_\C$ for its 
complexification. 

\par On the group level we let $G_\C$ be a simply connected 
Lie group with Lie algebra $\g_\C$. Write $A$, $A_\C$, $G$, $K$, $K_\C$, 
$N$ and $N_\C$ for the analytic subgroups of $G_\C$ with 
Lie algebra $\a$, $\a_\C$, $\g$, $\k$, $\k_\C$, $\n$ and $\n_\C$. 
The {\it Weyl group} of $\Sigma$ is given by  ${\Cal W}=N_K(\a)/ Z_K(\a)$.

\par For every $\alpha\in\Sigma$ we define a character 
$\xi_\alpha$ of $A_\C$ by 
$$\xi_\alpha(a)=a^\alpha\:=e^{\alpha(\log a)}\ .$$
{\it Regular} and {\it singular} 
elements in $A_\C$ are defined by 
$$A_{\C, \text{reg}}=\{ a\in A_\C \: 
(\forall \alpha\in \Sigma)\ \xi_\alpha^2(a)\neq 1\}\quad 
\hbox{and}\quad A_{\C, \text{sing}}=A_\C \bs A_{\C, \text{reg}}\ .$$

\subhead 1.2. Algebraic aspects and categorical quotients \endsubhead 
Write $X=G/K$ for the Riemannian symmetric space 
associated with $G$ and $K$. We denote by $X_\C=G_\C/ K_\C$ the 
complexification of $X$. Notice that $X_\C$ is an affine $K_\C$-variety.
As we are interested in the $K_\C$-action 
on $X_\C$ an algebraic view point is useful for our geometric studies.
\par The results collected below are not new, but sometimes 
we include short proofs for readability. We first recall the basic 
facts concerning affine varieties 
with reductive group action. 

\par Let $H$ be a complex reductive group and $V$ be an 
affine $H$-variety. Define {\it semisimple} and {\it regular} elements in 
$V$ by 

$$V_s=\{ v\in V\: Hv\ \hbox{is closed}\}$$
and 
$$V_r=\{ v\in V\: \dim Hv \ \hbox{is maximal}\}\ .$$
Further we set $V_{rs}=V_r\cap V_s$.
\par Denote by $\C[V]$ the algebra of regular (polynomial)
functions on $V$ and write $\C[V]^H$ for the $H$-invariants 
in $\C[V]$. According to Hilbert,  $\C[V]^H$ is finitely generated 
and so corresponds to an affine variety $V//H$, the {\it categorical quotient of $V$ by $H$}.  Let us 
denote by 
$$P\: V\to V//H$$
the map corresponding to the inclusion $\C[V]^H \into \C[V]$.
Recall from [6, II.3.2-3] the following properties of $P$.

\proclaim{Lemma 1.2.1} The quotient map $P\: V\to V//H $ has the following 
properties:
\roster 
\item Every fiber of $P$ contains a unique closed orbit. 
In particular $P$ is onto.
\item $(\forall v\in V_s)\ P^{-1}(P(v))=\{w\in V\:v\in \oline{ Hw}\}$.
\item $(\forall v\in  V_{rs}) \ P^{-1}(P(v))=Hv$.
\endroster
\qed 
\endproclaim

Let us now switch to our special situation of $V=X_\C$ and 
$H=K_\C$. We denote by $x_o=K_\C$ the base point in $X_\C$. 
Then, as a special case of general results of Matsuki [10],
we obtain that:

\proclaim{Lemma 1.2.2} The following assertions hold: 
\roster
\item $X_{\C,s}=K_\C A_\C.x_o $.
\item $X_{\C, rs} =K_\C A_{\C, \text{\rm reg}}.x_o$.
\endroster \qed 
\endproclaim

We continue with a discussion of the set $X_{\C,rs}$. 
Define a finite $2$-group by $F=A_\C \cap K_\C$. Notice that 
${\Cal W}$ normalizes $F$ and that both sets $A_{\C, \text{\rm reg}}$ and
$A_{\C, \text{\rm sing}}$
are invariant under multiplication by $F$.
Write $M_\C=Z_{K_\C}(A)$. 
We define a right action of ${\Cal W}$ on 
$K_\C/M_\C\times A_{\C, \text{\rm reg}}/F$ by $(kM_\C, aF).w=(kwM_\C, w^{-1} aw F)$. 
The resulting quotient will be denoted by $K_\C/M_\C\times_{\Cal W} A_{\C, \text{\rm reg}}/F$. Then it was essentially shown in 
[6, Lemma 2.1] that:

\proclaim{Proposition  1.2.3} The set $X_{\C, rs}$ is open and dense in 
$X_\C$. Moreover, the mapping 
$$\Phi\:  K_\C/M_\C\times_{\Cal W} A_{\C, \text{\rm reg}}/F\to X_{\C, rs}, \ \ 
[kM_\C, aF]\mapsto ka.x_o$$
is a $K_\C$-equivariant biholomorphism.
\qed \endproclaim 

\medskip In the remainder of this subsection we discuss
more specifically the structure of the quotient mapping 
$P\: X_\C \to X_\C//K_\C$.

\par Recall that the algebra $\C[X_\C]^{K_\C}$ of invariants is free 
in $n=\dim_\R \a$ generators, say 
$$\C[X_\C]^{K_\C}=\C[p_1, \ldots, p_n]\simeq\C [x_1, \ldots, x_n]\ .$$
Thus $X_\C// K_\C=\C^n$ and the projection mapping is 
given by 

$$P\: X_\C\to \C^n, \ \ z\mapsto (p_1(z), \ldots, p_n(z))\ .$$ 

Notice that $A_\C/F\to A_\C.x_o, \ aF\mapsto a.x_o$ is a 
biholomorphism . We denote by $P_A\: A_\C.x_o\to \C^n$ the restriction of 
the projection $P$ to $A_\C.x_o$. We summarize some (known) facts 
of the map $P_A$:

\proclaim{Proposition 1.2.4} The restriction $P_A\: A_\C.x_o\to \C^n$ 
of $P$ to $A_\C.x_o$ has the following properties: 
\roster
\item The map $P_A$ is holomorphic onto open finite and proper. 
\item For every $z\in A_\C.x_o$ the fiber of $P_A$ 
over $P_A(z)=P(z)$ is ${\Cal W}.z$. 
\item The restriction of $P_A$ to the regular elements
$$P_A\res_{A_{\C, \text{\rm reg}}.x_o}\: A_{\C, \text{\rm reg}}.x_o \to 
P(A_{\C, \text{\rm reg}}.x_o)$$ is a covering map. 
\endroster
\endproclaim 

\demo{Proof} As $P$ is a polynomial map, its restriction $P_A$ is
holomorphic. Next $P_A$ is onto because of Lemma 1.2.1(1)
and Lemma 1.2.2(1). We do not show that $P_A$ is proper 
but remark that it is a known fact which can be deduced 
from the fact that the polynomials  $p_1, \ldots, p_n$ are 
the elementary spherical functions. Standard finite dimensional 
representation theory then gives the properness of $P_A$.

\par By standard structure theory we have for $z, z'\in A_\C.x_o$ 
that 
$$K_\C z=K_\C z'\qquad \iff  \quad (\exists w\in {\Cal W}) \ w.z=z'\ .$$  
Thus Lemma 1.2.1(1) and Lemma 1.2.2(1) give us  that $P_A$
is finite as well as the second assertion 2. in the proposition. 
\par That $P_A$ is open now follows from the properness 
and 2. This completes the proof of 1. 
\par Finally, the third assertion 3. follows from 1. and 2. 
\qed 
\enddemo

If $M$ is a complex manifold, then we write ${\Cal O}(M)$ for the 
space of holomorphic functions on $M$. If in addition a group 
$L$ acts on $M$, then we denote by ${\Cal O}(M)^L$ the subspace 
of ${\Cal O}(M)$ which consists of $L$-invariant functions. 
Then we obtain from Proposition 1.2.4 the following 
(known) corollary:

\proclaim {Corollary 1.2.5} Let $M\subeq X_\C $ be a connected open 
$K_\C$-invariant subset. Then  $M_A=M\cap A_\C.x_o$ is open
and ${\Cal W}$-invariant in $A_\C.x_o$  and the restriction mapping 
$$\Res\: {\Cal O}(M)^{K_\C}\to {\Cal O}(M_A)^{\Cal W}, \ \ 
f\mapsto f\res_{M_A}$$
is an isomorphism. 
\endproclaim  

\demo{Proof} Let us first show that $\Res $ is injective. Let 
$f\in {\Cal O}(M)^{K_\C}$ such that $f\res_{M_A}=0$. The $K_\C$-invariance 
of $f$ implies that 
$f$ vanishes on $K_\C M_A=M\cap X_{\C,s}$. As $X_{\C,s}\cap M$
is dense in $M$ by Proposition 1.2.3, it follows that 
$f\equiv 0$. Hence $\Res$ is injective. 

\par  To show that $\Res $ is onto, let $F\in {\Cal O}(M_A)$. We have to 
show that there exists an $f\in {\Cal O}(M)^{K_\C}$ such that 
$f\res_{M_A}=F$. Let us define a function on $M$ by  
$$(\forall z\in M)\qquad f(z)=F(P_A^{-1}(P(z)))\ .$$
This map $f$ is well defined as every fiber of $P_A$ is 
a ${\Cal W}$-orbit (cf. \ Proposition 1.2.4(2))
and $F$ is invariant under ${\Cal W}$.  The properness of $P_A$ 
gives us that $f$ is continuous. It is also clear that 
$f\res_{M_A}=F$ and that $F$ is $K_\C$-invariant. It remains 
to see that $f$ is holomorphic. {} From Proposition 1.2.4(3) 
we deduce that $f$ is holomorphic on $X_{\C, rs}\cap M$. As 
$f$ is also continuous, the Riemann removable singularity theorem 
applies and yields the holomorphicity of $f$. 
\qed 
\enddemo

\subhead 1.3. Complex crowns \endsubhead 
Define a ${\Cal W}$-invariant, convex and relatively compact subset of $\a$ by 
$$\Omega=\{ Y\in \a\: (\forall \alpha\in \Sigma)\ |\alpha(Y)|<{\pi\over 2}\}
\ .$$ 
Following [1] we define the {\it complex crown} $\Xi$ of $X$ by 

$$\Xi=G\exp(i\Omega).x_o\subeq X_\C\ .$$
It was shown in [1] that $\Xi$ is an open $G$-invariant neighborhood of $X$ 
with proper $G$-action. Notice that $\Xi$ is independent of the choice of $\a$ hence canonically 
attached to $X$. Recently it was shown by many authors that $\Xi$ is actually Stein. 

\par For a subset $\omega\subeq 2\Omega$ we define the set 
$$T(\omega)=A\exp(i\omega)$$ 
and notice that $T_\omega$ is homeomorphic to $\a+i\omega$
via the exponential mapping $\exp\res_{\a_\C}\: \a_\C \onto A_\C$. 
Moreover, as $T(2\Omega)\cap F=\{1\}$ it follows that 
the mapping $T(2\Omega)\to X_\C , a\mapsto a.x_o$ is an injection. 

\par Another important feature of the crown is that 

$$\Xi\subeq N_\C T(\Omega).x_o \leqno(1.3.2)$$
(cf.\ [2]). 
In particular we have a well defined 
holomorphic mapping 
$$a\: \Xi\to T(\Omega), \ \ nu.x_o\mapsto u \qquad 
(n\in N_\C, u\in T(\Omega))\ .$$
Further we have the complex convexity theorem [2] which states that 
$$(\forall Y\in \Omega)\qquad \Im \log a(G\exp(iY).x_o)
\subeq \conv({\Cal W}Y)\ .\leqno(1.3.3)$$
Here $\conv(\cdot)$ denotes the convex hull of $(\cdot)$. 

\medskip Set $\a^+=\{ Y\in \a\: (\forall \alpha\in \Sigma^+)\ \alpha(Y)>0\}$
and $A^+=\exp(\a^+)$. Later we will need also the 
following subsets of $\Omega$:

$$\Omega^+=\Omega\cap \a^+\quad \hbox{and} \quad \Omega'={\Cal W}\Omega^+\ .$$

\par It is useful to introduce some further terminology 
and define for very ${\Cal W}$-invariant convex subset 
$\omega\subeq 2\Omega$ the $K_\C$-invariant set 
$$X_{\C, \omega}=K_\C T(\omega).x_o$$
in $X_\C$. Notice that even for $\omega$ open the set  
$X_{\C,\omega}$ is not open in $X_\C$. To overcome this difficulty 
we define 
$$\hat X_{\C, \omega}=P^{-1}(P(T(\omega).x_o))\ .$$
{} From the definition it follows immediately that 
$\hat X_{\C, \omega}$ is a $K_\C$-invariant subset of $X_\C$ 
containing $X_{\C, \omega}$. Moreover it follows from 
Proposition 1.2.4(1) that  $\hat X_{\C, \omega}$ is open if
$\omega$ is open.

\head 2. The convexity theorem\endhead 

For $Y\in \a$ let us denote by $\co(Y)$ the convex hull of the Weyl group 
orbit ${\Cal W} Y$.  We start with stating the main result.

\proclaim {Theorem 2.1} Let $Y\in \Omega$. Then 
$$G\exp(iY).x_o\bigcap X_{\C, s}\subeq X_{\C, \co(Y)}=K_\C T(\co(Y)).x_o\ .\leqno(2.1)$$
In particular,  
$$G\exp(iY).x_o\subeq \hat X_{\C, \co(Y)}\ .\leqno(2.2)$$
\endproclaim 

Before we start proving Theorem 2.1 let us mention an interesting 
consequence:

\proclaim {Corollary 2.2} The following inclusions hold: 
\roster 
\item $\Xi\bigcap X_{\C, s}\subeq X_{\C, \Omega}.$ 
\item $\Xi\subeq \hat X_{\C, \Omega}$. 
\endroster\qed 
\endproclaim 

It seems to us that giving a direct geometric proof of Theorem 2.1 is 
rather difficult. Instead we will proceed indirectly. We 
will use spherical functions on $X$. These functions have an analytic continuation 
to a holomorphic function on $\Xi$ as well as to a 
$\hat X_{\C, 2\Omega}$  (cf.\ [8]). Below we will 
give pointwise estimates of the analytically extended spherical functions
which are uniform in the spectral parameter. From that we will deduce 
Theorem 2.1.

\subhead 2.1. Spherical functions\endsubhead 
We begin with recalling Harish-Chandra's integral representation 
of spherical functions. For $\alpha\in \Sigma$ we set 
$m_\alpha=\dim \g^\alpha$ and define $\rho={1\over 2}\sum_{\alpha\in 
\Sigma^+} m_\alpha\alpha\in \a^*$. 
\par If $\lambda \in \a_\C^*$, then one 
defines the {\it  spherical function on $X$ 
with parameter $\lambda$}  to be the $K$-invariant eigenfunction 

$$\phi_\lambda(gK)=\int_K a(kg)^{\rho-\lambda} \ dk \qquad (g\in G)\ .$$
Here $dk$ is the normalized Haar measure on the compact group $K$. 
It follows from (1.3.2) that $\phi_\lambda$ has a holomorphic 
extension to $\Xi$ which we will denote by $\tilde \phi_\lambda$ (cf.\ [8]). 
\par On the other hand $\phi_\lambda\res_A$ extends to 
a ${\Cal W}$-invariant holomorphic function on $T(2\Omega)$ [8].  
Thus Corollary 1.2.5 applied to $M=\hat X_{\C, 2\Omega}$ implies:

\proclaim {Proposition 2.1.1} Let $\lambda\in \a_\C^*$. Then $\phi_\lambda$
uniquely extends to a $K_\C$-invariant holomorphic function 
on $\hat X_{\C, 2\Omega}$.\qed \endproclaim

We shall prove point-wise estimates for analytically continued spherical functions
which are uniform in the spectral paramter. The upper bound is 
easy to obtain:

\proclaim {Lemma 2.1.2} Let $Q\subeq \Xi$ be a compact subset. 
Then there exists a constant $C>0$ such that for all $g\exp(iY).x_o\in Q$
with $g\in G, Y\in \Omega$ one has 
$$(\forall \lambda\in i\a^*)\qquad 
|\tilde\phi_\lambda(g\exp(iY).x_o)|\leq C \cdot \sup_{w\in {\Cal W}} e^{i\lambda(wY)}\ .
$$
\endproclaim 

\demo{Proof} Fix $g\exp(iY).x_o\in Q$ and recall the definition 

$$\tilde\phi_\lambda (g\exp(iY).x_o)=\int_K a (kg\exp(iY))^{\rho-\lambda}\ dk \ .$$
Thus (1.3.3) implies that  
$$|\tilde\phi_\lambda (g\exp(iY).x_o)|\leq \int_K 
\left|a (kg\exp(iY))^{\rho-\lambda}\right |\ dk \leq 
C \cdot \sup_{w\in {\Cal W}} e^{i\lambda(wY)}$$
with 
$$C=\sup_{z\in Q} \int_K \left|a (kz)^\rho\right|\ dk<\infty\ .$$ 
\qed 
\enddemo

Now we switch to  the lower estimates. 
For that we will use Harish-Chandra's power series expansion of spherical function (cf.\ [4]). 

\par We briefly recall the standard notation. 
Define a semigroup $\Lambda=\N_0[\Sigma^+]\subeq \a^*$. 
As usual we denote by $\c(\lambda)$ the Harish-Chandra $c$-function. 
This is a meromorphic function on $\a_\C^*$ 
with no poles on the imaginary axis $i\a^*$; it can be explicitely
computed (Gindikin-Karpelevic formula, see [4, Ch.\ IV, Th.\ 6.14]).   

\par For  $\sigma\in \a_\C^*$ and $a\in T(\Omega)$ we write $a^\sigma=e^{\sigma(\log a)}$. 

\par Fix $\lambda\in i\a^*$. Then for $a\in A^+$ Harish Chandra's 
expansion of $\phi_\lambda$ is given as follows
(cf.\ [4, Ch.\ IV, Th. 5.5]): 
$$\phi_\lambda(a)=\sum_{w\in {\Cal W}} \c(w\lambda) \cdot a^{w\lambda-\rho}
\cdot \Psi_{w\lambda}(a)\ , \leqno(2.1.1)$$
where 
$$\Psi_{w\lambda}(a)=  \sum_{\mu\in \Lambda}
\Gamma_\mu(w\lambda) a^{-\mu}\ .$$

For $\sigma\in i\a^*$ the coefficients $\Gamma_\mu(\sigma)$ 
are defined recursively by $\Gamma_0(\sigma)=1$ and 

$$\Gamma_\mu(\sigma)={2\over \la \mu, \mu\ra +2\la \mu, \sigma\ra}
\sum_{\alpha\in \Sigma^+} m_\alpha \sum_{k\geq 1}
\Gamma_{\mu-2k\alpha}(\sigma) (\la \mu+\rho-2k\alpha, \alpha\ra +
\la \alpha, \sigma\ra )\ .\leqno (2.1.2)$$

The power series $\Psi_{w\lambda}$ are locally 
abolutely convergent on $A^+$. One can deduce this 
from the following fact ([4, Ch.\ IV, Lemma 5.6]): 
for all $H\in \a^+$ there is a constant $C_H>0$ such that 

$$(\forall \mu\in \Lambda)(\forall \lambda\in i\a^*)\qquad |\Gamma_\mu(\lambda)|\leq C_H\cdot 
e^{\mu(H)}\ .\leqno (2.1.3)$$
It is clear that (2.1.1) is also valid for all $a\in A^+\exp(i\Omega)$. 

\proclaim {Lemma 2.1.3} Fix $\lambda_0\in i\a^*$. 
Then the following assertions hold:
\roster 
\item For all $\mu\in \Lambda$ the limit 
$$\Gamma_\mu[\lambda_0]\:=\lim_{t\to \infty}\Gamma_\mu(t\lambda_0)$$ 
exists. 
\item The power series 
$$\Psi_{[\lambda_0]}(a)=\sum_{\mu\in\Lambda} \Gamma_\mu([\lambda_0])
a^{-\mu}\ \qquad (a\in A^+\exp(i\Omega)) $$ 
is locally absolutely convergent on $ A^+\exp(i\Omega)$. 
In particular, $\Psi_{[\lambda_0]}$ defines 
a non-zero holomorphic function on $ A^+\exp(i\Omega)$. 
Moreover 
$$\Psi_{t\lambda_0}\to \Psi_{[\lambda_0]}$$
locally uniformly on $A^+\exp(i\Omega)$ for $t\to \infty$. \endroster
\endproclaim 

\demo{Proof} 1. Let $\alpha_1,\ldots, \alpha_n$ be a basis of $\Sigma^+$. 
If $\mu=m_1\alpha_1+\ldots+m_n\alpha_n\in \Lambda $, then define 
its degree by $\deg \mu=m_1+\ldots+m_n$. We use the recurrence 
relation (2.1.3) and induction on $\deg\mu$  to deduce that 
for every $\mu\in \Lambda$ the function 
$$F_\mu\: \R \to \C , \  \ t\mapsto \Gamma_\mu(t\lambda_0)$$
is rational with no poles. Moreover (2.1.3) implies that $F_\mu$ is 
bounded. Thus if $F_\mu={p_\mu\over q_\mu}$ with polynomials
$p_\mu$ and $q_\mu$ we obtain that $\deg p_\mu\leq \deg q_\mu$. 
From that the assertion follows. 

\par\nin 2. This follows from 1. and the estimate (2.1.3).\qed 
\enddemo

Let us call an element $\lambda\in \a_\C^*$ {\it regular}
if $\la \lambda, \alpha \ra \neq 0$ for all $\alpha\in \Sigma$.

\proclaim{Proposition 2.1.4} Fix $\lambda_0\in i\a^*$ regular. 
Let $H\in \a^+$ and $Y\in \Omega'$ and set $a=\exp(H+iY)$. 
Then there exists a constant $C_{a,\lambda_0}\geq 0$ such that 
for all $t>>0$ one has 
$$|\tilde\phi_{t\lambda_0}(a.x_o)|\geq {C_{a,\lambda_0}\over t^p}\cdot 
\sup_{w\in {\Cal W}} e^{it\lambda_0(wY)}, $$
where $p={1\over 2}\dim N$. Moreover 
$C_{a,\lambda_0}>0$ for all elements $a$ in a dense open subset 
of $A^+\exp(i\Omega')$. \endproclaim

\demo{Proof} A basic result of Harish-Chandra asserts 
that $\phi_\lambda=\phi_{w\lambda}$ for all $\lambda\in \a_\C^*$ 
and $w\in {\Cal W}$. Thus by the regularity of $\lambda_0$ and 
$Y$ we may assume that $\lambda_0(iY)> \lambda_0(iwY)$ for all 
$w\in {\Cal W}$, $w\neq \1$. 
Therefore it suffices to show that 
$$(\forall t>>0)\qquad |\tilde\phi_{t\lambda_0}(a.x_o)|\geq {C_{a,\lambda_0}\over t^p}\cdot 
e^{it\lambda_0(Y)}\leqno(2.1.4)$$
for a constant $C_{a,\lambda_0}\geq 0$ which is actually positive 
for $a$ ranging in a dense open subset of $A^+\exp(i\Omega')$.  

\par {} From (2.1.1) we obtain that: 
$$\tilde\phi_{t\lambda_0}(a.x_o)=a^\rho\cdot \sum_{w\in {\Cal W}} \c (tw\lambda_0)
\cdot e^{it\lambda_0(wY)}\cdot e^{t\lambda_0(wH)} \cdot
\Psi_{tw\lambda_0}(a)\ . \leqno(2.1.5)$$
Contemplating on (2.1.5) we notice the following facts:
\medskip

\item{$\bullet$} $|e^{t\lambda_0(wH)}|=1$ for all $t\in \R$ and 
$w\in {\Cal W}$
as $\lambda_0$ is imaginary. 
\item{$\bullet$} $|a^\rho|>0$ is independent 
of $\lambda$. 
\item{$\bullet$} $|\c(tw\lambda_0)|=|\c(t\lambda_0)|$ for all 
$w\in {\Cal W}$ (Maa\3-Selberg relations, cf.\ [4, p.\ 451, (16)]). 
\item{$\bullet$} There is a constant $C>0$ such that 
$|\c(t\lambda_0)|\geq {C\over t^p}$
for $t>>0$ (cf.\ [4, Ch.\ IV, Prop.\ 7.2]). 
\item {$\bullet$} $\Psi_{tw\lambda_0}(a)$ is uniformly  bounded in 
$t$ for all $w\in {\Cal W}$ by Lemma 2.1.2 (2). 
\medskip 
Combining these facts with (2.1.5) we obtain (2.1.4) for some 
$C_{a,\lambda_0}\geq 0$. To be more precise we have 
$C_{a,\lambda_0}>0$ provided  that $\Psi_{[\lambda_0]}(a)\neq 0$ 
(cf.\ Lemma 2.1.3(2)),  and this 
finishes the proof.  
\qed 
\enddemo

\subhead 2.2. Proof of Theorem 2.1 \endsubhead 
During the proof we will frequently use the fact that 
the mapping 
$$G/M\times \Omega^+\to \Xi, \ \ (gM, Y)\mapsto g\exp(iY).x_o\leqno (2.2.1)$$
is a diffeomorphism with open and dense image (cf.\ [9, Cor.\ 4.2]). 
\par Due to the statement of the theorem we may assume that 
$Y\in \Omega^+$. We now proceed by contradiction. 
If the statement were false, there would exist
a $g\in G$ such that 
$g\exp(iY).x_o\in X_{\C,s}$ but 
$g\exp(iY).x_o\not\in K_\C T(\co(Y)).x_o$. 
Let $g(t)$, $t\in [0,1]$,  be a continuous path in $G$ such that $g(0)=\1$ and 
$g(1)=g$. Recall the quotient mapping $P\: X_\C \to \C^n$ from 
Subsection 1.2. Define a continuous path in $\C^n$ by $\gamma(t)=P(g(t)\exp(iY).x_o)$. 
Then $\gamma(0)=P(\exp(iY).x_o)\in P(T(\co(Y)).x_o)$. 
As $P_A=P\res_{A_\C.x_o}$
is proper (cf.\ Proposition 1.2.4(1)), we conclude 
that $P(T(\co(Y)).x_o)$ is closed. Thus there exists a $0\leq t_0< 1$
which is maximal with respect to the property that 
$\gamma(t)\in P(T(\co(Y)).x_o)$ for all $0\leq t\leq t_0$. 
Again it follows from Proposition 1.2.4(1)  that $P(T(\Omega).x_o)$ 
is an open neighborhood of $\gamma(t_0)$. Thus we find a $t_0<t_1\leq 1$ 
such that $\gamma(t)\in P(T(\Omega).x_o)$ for all $0\leq t\leq t_1$. 
Now pick  $t_0<t_2\leq t_1$ such that $\gamma(t_2)=P(a.x_o)$
for some $a\in T(\Omega)\bs T(\co(Y))$. It is no loss 
of generality to assume now that $g=g(t_2)$ and reparametrizing $g(t)$ we  
may also assume that $t_2=1$. 

\par To sum up, if the statement were false we find a $g\in G$
and a  continuous path 
$g(t)$, $t\in [0,1]$, in $G$  with $g(0)=\1$ and $g(1)=g$ 
such that 
$g(t)\exp(iY).x_o\in \hat X_{\C, \Omega}$ for all $t\in [0,1]$
and $P(g\exp(iY).x_o)=P(a.x_o)$ for some $a\in T(\Omega)\bs 
T(\co(Y))$. In particular, it follows from 
Proposition 2.1.1 that 

$$(\forall \lambda\in \a_\C^*)\qquad 
\tilde \phi_\lambda(g\exp(iY).x_o)=\tilde\phi_\lambda(a.x_o)\ . 
\leqno (2.2.2)$$

\par Using (2.2.1) we  can vary $g$ and $Y$ 
slighty to obtain $a\in A^+\exp(i\Omega')$, say 
$a=\exp(H+iZ)$
for some $H\in \a^+$ and $Z\in \Omega'$. 

\par As $Z\in \Omega \bs \co(Y)$,  we can use the theorem of 
Hahn-Banach to find a $\lambda_0\in i\a^*$ regular such that 
$|\lambda_0(Z)|>|\lambda_0(wY)|$ for all $w\in {\Cal W}$. 
\par Lemma 2.1.2 implies 
$$(\forall t>0)\qquad |\tilde\phi_{t\lambda_0}(g\exp(iY).x_o)|\leq C \cdot 
\sup_{w\in {\Cal W}} e^{ti\lambda_0(wY)}\leqno 
(2.2.3)$$
for a constant $C>0$. 
\par On the other hand (2.2.2) 
combined with Proposition 2.1.4
gives us for all $t>>0$ that 

$$  |\tilde\phi_{t\lambda_0}(g\exp(iY).x_o)|= |\tilde\phi_{t\lambda_0}(\exp(H+iZ).x_o)|
\geq {C'\over t^p} \cdot \sup_{w\in W}  e^{ti\lambda_0(wZ)}\leqno (2.2.4)$$
for a constant $C'\geq 0$. Moreover varying $g$ and $Y$ slightly if necessary 
we may assume that $C'>0$ (cf.\ Proposition 2.1.4). Then (2.2.4) clearly 
contradicts (2.2.3),  concluding the proof of the theorem.\qed

\subhead 2.3. Applications\endsubhead  
Define a subset of $\a_\C^*$ by 
$$\Pi=\{ \lambda\in \a_\C^*\: \phi_\lambda \ \hbox{is positive definite}\}\ .$$
Recall that $i\a^*\subeq \Pi$ and that for all $\lambda\in\Pi$ 
the spherical function $\phi_\lambda$ is bounded on $X$. 

\par Using Theorem 2.1 and the same reasoning as in [7] we now obtain 
that:

\proclaim {Theorem 2.3.1} Let $\lambda\in \Pi$. Then the analytically 
continued spherical function $\tilde\phi_\lambda$ is bounded 
on $\Xi$.\qed
\endproclaim 

\demo{Proof} For $\lambda\in \Pi$ it was established 
in [8] that $\tilde\phi_\lambda$ is bounded on $X_{\C,\Omega}$. 
Now the assertion follows from Corollary 2.2. \qed 
\enddemo

Another application of Theorem 2.1 is related to the heat
kernel which was suggested to us by M. Stenzel. 
In the sequel $t>0$ denotes a positive number. 
We denote by $\rho_t(x)$, $x\in X$,   the heat kernel on $X$. 

Using the spectral resolution of $\rho_t$ it was shown in [9] that 
$\rho_t$ has an analytic continuation to a holomorphic function 
on $\Xi$ which we will denote by $\tilde\rho_t$. 
Then it follows readily from Theorem 2.1, Lemma 2.1.2 
and the spectral resolution of $\rho_t$ that:

\proclaim {Theorem 2.3.2} Fix $t>0$ and let $\rho_t$ be the heat 
kernel on $X$.  Then the analytically 
continued heat kernel $\tilde\rho_t$ is bounded on $\Xi$. 
\qed 
\endproclaim

\head 3. Examples:  Symmetric matrices\endhead
In this section we will illustrate Theorem 2.1 for the choices 
of $G=\Sl(n,\R)$ and $G=\Sl(n,\C)$. Then $X=G/K$ has 
a canonical realization in real/complex symmetric matrices.  
It might be interesting to observe that even in these classical cases 
Theorem 2.1 seems to be new and likely forfeits a proof using only techniques 
from elementary linear algebra.

\subhead 3.1. Real symmetric matrices \endsubhead 
Let now $G=\Sl(n,\R)$ and $K=\SO(n)$. Then $G_\C=\Sl(n,\C)$
and $K_\C=\SO(n,\C)$. The  mapping 
$gK_\C\mapsto gg^t$ identifies $X_\C$ with the
complex symmetric matrices of determinant one
$$X_\C=\{ Z\in M(n,\C)\: Z=Z^t, \det Z=1\}\ .$$
Within this identification $X$ is the set 
of all real positive definite 
symmetric matrices of determinant one

$$X=\{ Y\in M(n,\R)\: Y=Y^t, Y>>0, \det Y=1\}\ .$$

\par We choose $A$ to be the positive diagonal matrices in $G=\Sl(n,\R)$. 
Then ${\Cal W}$ becomes the symmetric 
group ${\frak S}_n$ in $n$ letters permuting the diagonal 
entries.  The set 
$\Omega$ is given by 
$$\Omega=\{ \diag(y_1, \ldots, y_n)\in M(n,\R)\: 
\sum_{i=1}^n y_i=0, \ |y_i-y_j|<{\pi\over 2}\}\ .$$
Accordingly we have 
$$T(\Omega)=\{a=\diag(a_1,\ldots,a_n)\in \Sl(n,\C)\: 
(\forall 1\leq i,j\leq n)\quad \Re \left( {a_i\over a_j}\right)>0\}\ .$$
Let us write $D_\R$ for the diagonal matrices in $X$ and 
$D_\C$ for the diagonal matrices in $X_\C$. Notice that 
the mapping $A_\C\to D_\C, a\mapsto a.x_o=a^2$ is the square mapping. 
If $D(\Omega)$ denotes the isomorphic image of $T(\Omega)$ in $X_\C$, 
then we have 
$$D(\Omega)=\{a=\diag(a_1,\ldots,a_n)\in X_\C\: 
(\forall 1\leq i,j\leq n)\quad  {a_i\over a_j}\in \C_+\}\ , $$
where $\C_+=\C\bs ]-\infty,0]$. 

Notice that 
$$\Xi=\bigcup_{g\in G} gD(\Omega)g^t\ .$$
\par Now Theorem 2.1 in this special situation
reads as follows:

\proclaim{Proposition 3.1.1} Let $Y=\diag(y_1, \ldots, y_n)$ be a real 
diagonal matrix with zero trace and $|y_i-y_j|<\pi$ for all $1\leq i,j\leq n$. Then 
$$\co(Y)=\conv\{ \diag (y_{\sigma(1)}, \ldots, y_{\sigma(n)})\:
\sigma\in {\frak S}_n\}$$
and 
$$\left[\bigcup_{g\in \Sl(n,\R)} ge^{iY} g^t\right]\bigcap
\left[\bigcup_{k\in \SO(n,\C)} k D_\C k^{-1}\right]
\subeq \bigcup_{k\in \SO(n,\C)} k D_\R e^{i\co(Y)} k^{-1}\ .$$
\qed \endproclaim

\subhead 3.2. Complex symmetric matrices \endsubhead 
Let us now consider the case of $G=\Sl(n,\C)$ and $K=\SU(n)$. 
Then we have canonical biholomorphisms 
$K_\C \simeq G$ and $X_\C\simeq G$. Within this identifications
the actions of $G$ and $K_\C$ on $X_\C$ are given by 

$$G\times X_\C\to X_\C, \ \ (g,z)\mapsto gzg^*\ , $$
$$K_\C\times X_\C \to X_\C, \ \ (k,z)\mapsto kzk^{-1}$$   
(cf.\ [3, Sect. 3]). As before we choose $A<G$ to be the 
group of diagonal matrices with positive entries and write 
$D(\Omega)$ for the faithful image of $T(\Omega)$ in $X_\C=G$. 
Then 

$$D(\Omega)=\{a=\diag(a_1,\ldots,a_n)\in X_\C\: 
(\forall 1\leq i,j\leq n)\quad  {a_i\over a_j}\in \C_+\} .$$
In particular we obtain 

$$\Xi=\bigcup_{g\in \Sl(n,\C)} gD(\Omega)g^*\ .$$

\par With the notation introduced above  Theorem 2.1 now reads as follows:

\proclaim{Proposition 3.2.1} Let $Y=\diag(y_1, \ldots, y_n)$ be a real 
diagonal matrix with zero trace and $|y_i-y_j|<\pi$ for all $1\leq i,j\leq n$. Then 
$$\left[\bigcup_{g\in \Sl(n,\C)} ge^{iY} g^*\right]\bigcap
\left[\bigcup_{k\in \Sl(n,\C)} k D_\C k^{-1}\right]
\subeq \bigcup_{k\in \Sl(n,\C)} k D_\R e^{i\co(Y)} k^{-1}\ , $$
where $\co(Y)=\conv\{ \diag (y_{\sigma(1)}, \ldots, y_{\sigma(n)})\:
\sigma\in {\frak S}_n\}$. 
\qed \endproclaim

\refstyle{C}

\enddocument